\documentclass[trsc, nonblindrev]{informs3}

\OneAndAHalfSpacedXI 

\usepackage{preamble}

\usepackage{natbib}
 \bibpunct[, ]{(}{)}{,}{a}{}{,}%
 \def\bibfont{\small}%
\TheoremsNumberedThrough     

\EquationsNumberedThrough    

\MANUSCRIPTNO{TS-2024-0604}

\begin{document}


\RUNAUTHOR{Wouda, Romeijnders, and Ursavas}
\RUNTITLE{Feasibility cuts for chance-constrained multicommodity network design}
\TITLE{On feasibility cuts for chance-constrained multicommodity network design problems}

\ARTICLEAUTHORS{%
\AUTHOR{Niels A. Wouda, Ward Romeijnders, Evrim Ursavas}
\AFF{Department of Operations, University of Groningen, \EMAIL{\{n.a.wouda, w.romeijnders, e.ursavas\}@rug.nl} \URL{}}
} 

\ABSTRACT{%
\textbf{Problem definition}:
We study efficient exact solution approaches to solve chance-constrained multicommodity network design problems under demand uncertainty, an important class of network design problems.
The chance constraint requires us to construct a network that meets future commodity demand sufficiently often, which makes the problem challenging to solve.
\textbf{Methodology/results}:
We develop a solution approach based on Benders' decomposition, and accelerate the approach with valid inequalities and cut strengthening.
We particularly investigate the effects of different subproblem formulations on the strength of the resulting feasibility cuts.
We propose a new formulation that we term \emph{FlowMIS}, and investigate its properties.
Additionally, we numerically show that FlowMIS outperforms standard formulations: in our complete solution approach with all enhancements enabled, FlowMIS solves 67 out of 120 solved instances the fastest, with an average speed-up of 2.0\texttimes~over a basic formulation.
\textbf{Implications}:
FlowMIS generates strong feasibility cuts tailored to subproblems with a network flow structure.
This results in reduced solution times for existing decomposition-based algorithms in the context of network design, and the ability to solve larger problems.
}%


\maketitle

%

\section{Introduction}
\label{sec:introduction}

Network design models see many uses in industry, particularly in transportation and logistics~\citep{cordeau_survey_1998,crainic_service_2000,crainic_taxonomy_2022}, in telecommunications~\citep{wong_telecommunications_2021}, and in the energy sector~\citep{huang_integrated_2014,uster_biomass_2018,fragkos_decomposition_2021}.
These models often incorporate some form of strategic decision-making, involving the construction of new infrastructure such as the placement of facilities, cell towers or power lines, aimed to satisfy future demand for the network infrastructure.
Such future demand is typically uncertain.

In this paper we consider multicommodity capacitated fixed-charge network design problems with uncertain demand.
Our aim is to construct a network such that sufficient amounts of each commodity can flow through the network to meet the uncertain demand.
Rather than requiring demand to be met for all realisations, as in~\cite{rahmaniani_accelerating_2018}, we instead impose a service level constraint in the form of a chance constraint, which instead states that demand must be satisfied for sufficiently many different realisations.
Thus, it is permissible that there is a (typically small) subset of realisations for which demand cannot be satisfied, as long as the chance constraint is respected.
Service level guarantees are commonly used in industry to hedge against uncertainty, but models incorporating them are notoriously difficult to solve.

Since chance-constrained multicommodity capacitated fixed-charge network design problems are difficult to solve exactly using standard methods, we develop a tailored solution method.
We assume that the uncertainty can be captured according to a finite number of scenarios.
Our solution method is based on Benders' decomposition: we iteratively solve a master problem to determine a solution to the network design problem, and then determine whether the solution is feasible across sufficiently many scenarios by solving scenario subproblems that have a network flow structure.
To obtain good solutions already in early iterations, we strengthen the master problem with valid inequalities.
When insufficient scenarios are feasible, we derive feasibility cuts from the infeasible scenario subproblems, and add them to the master problem.
We show that different formulations of these scenario subproblems result in different feasibility cuts, and study their properties.
Since our subproblems have a network flow structure, it turns out that these properties relate closely to cuts in the network.
We propose a novel formulation, \textit{FlowMIS}, that generates strong feasibility cuts tailored to subproblems with a network flow structure.
These cuts significantly improve over formulations known from the literature.

Network design models under uncertainty have received increased attention in recent years.
\cite{keyvanshokooh_hybrid_2016} present enhancements of a Benders' decomposition algorithm in the form of valid inequalities, Pareto-optimal cuts in the manner of~\cite{magnanti_network_1984}, and a scenario reduction technique.
They apply their algorithm to a supply chain network design problem under demand and transportation uncertainties.
\cite{crainic_partial_2021} formalise the notion of \textit{partial Benders' decomposition}, whose core idea is to retain some scenario information in the master problem, rather than decompose all scenarios.
They consider both scenario retention, where existing scenarios are retained inside the master problem, and scenario creation, where artificial valid scenarios are generated and included inside the master problem.
They apply their partial decomposition algorithm to a stochastic multicommodity capacitated fixed-charge network design problem. 
\cite{rahmaniani_accelerating_2018} summarises these and many other enhancements in the context of a stochastic multicommodity fixed-charge network design problem, where all demand must be satisfied at minimal expected cost.
Most recently,~\cite{rahmaniani_asynchronous_2024} present an asynchronous parallel decomposition algorithm for stochastic multicommodity network design problems that effectively scales to multiple processor cores.

The literature on Benders' decomposition is reviewed in~\cite{rahmaniani_benders_2017}, and a large number of enhancements identified there are implemented in the context of stochastic network design in~\cite{rahmaniani_accelerating_2018}.
For our paper, we are particularly interested in two such enhancements: valid inequalities and cut strengthening.
Valid inequalities are often used to strengthen the master problem~\citep{belieres_benders_2020,fragkos_decomposition_2021}.
These additional constraints may obviate the need for feasibility cuts~\citep{peng_probabilistic_2020,crainic_partial_2021}, but that is not the case in our setting where we have subproblems particularly to test for feasibility.
Several studies also use Benders' decomposition to solve problems with pure feasibility subproblems~\citep{de_kruijff_integer_2018,fu_research_2019}.
An enhancement of Benders' feasibility cuts tailored to general chance-constrained mixed-integer problems with finite support is discussed in~\cite{luedtke_branch-and-cut_2014}, while~\cite{costa_benders_2009} present an enhancement specific to multicommodity network design problems.
Enhancements from the general branch-and-cut literature often apply to Benders' feasibility cuts as well.
As a particularly relevant example, reformulating the cut-generating subproblem to strengthen cuts is applied in~\cite{balas_mixed_1996} and~\cite{balas_modified_1997} to mixed-binary mathematical programs to be solved by branch-and-cut, and extended to Benders' feasibility cuts in~\cite{fischetti_note_2010}.
These subproblem reformulations apply to general problems, but stronger results can be obtained by tailoring the reformulation to the specific problem, as we will show.

We evaluate our proposed solution approach via a numerical study on 126 benchmark instances.
The results show that FlowMIS outperforms other formulations, solving the highest number of instances (120) when incorporating the valid inequalities and cut strengthening.
In comparison, other formulations all solve instances about as fast as a basic implementation of the decomposition algorithm, which FlowMIS outperforms by 2.0 times on average.
The experiments establish the efficiency of FlowMIS, particularly when augmented with additional enhancements.

Summarising, this paper makes the following contributions:
\begin{itemize}
    \item We develop a tailored solution method based on Benders' decomposition for chance-constrained multicommodity network design problems with uncertain commodity demand, incorporating valid inequalities and cut strengthening.
    \item We propose a novel subproblem formulation, FlowMIS, for generating strong feasibility cuts tailored to subproblems with a network flow structure.
    \item We prove properties of different scenario subproblem formulations, corresponding to our and various other formulations from the literature.
    \item Through a numerical study, we demonstrate the effectiveness of FlowMIS in solving a variety of benchmark instances.
\end{itemize}

The rest of this paper is structured as follows.
In~\cref{sec:problem_definition} we formally present our chance-constrained multicommodity capacitated fixed-charge network design problem.
In~\cref{sec:solution_approach} we develop our solution approach based on Benders' decomposition.
In~\cref{sec:numerical_experiments} we show this solution approach to be effective via a numerical study on $126$ benchmark instances.
Finally,~\cref{sec:conclusion} concludes the paper.

\vspace{1em}
\emph{Notation:}
Throughout this paper we use the following convenient notational shorthands.
First, we let $[N] = \{1, \ldots, N\}$ denote the set of natural numbers up to and including $N$.
Second, we let $\B = \{0, 1\}$ denote the set of binary digits.
Third, we use $\1{p(\cdot)}$ as an indicator variable that is $1$ when predicate $p(\cdot)$ is true, and $0$ otherwise.

\section{Problem definition}
\label{sec:problem_definition}

Following the notation of~\cite{chouman_commodity_2017} and~\cite{rahmaniani_accelerating_2018}, we consider a chance-constrained multicommodity fixed-charge network design problem on a directed graph $G = (N, A)$, where $N$ is the set of nodes and $A \subseteq N \times N$ the set of arcs.
There is a set of commodities $K$.
For each commodity $k \in K$ an uncertain amount of demand $d^k(\omega) \ge 0$ needs to be routed through the graph from an origin $O(k) \in N$ to a destination $D(k) \in N$.
The uncertain demand $d^k(\omega)$ depends on a random vector $\omega$ with known support $\Omega$.
To route these commodities through the graph, we first have to construct a network $y = (y)_{ij} \in \B^{|A|}$ for each arc $(i, j) \in A$ in the graph.
Each arc $(i, j) \in A$ is constructed at fixed cost $f_{ij} \ge 0$, and, when $y_{ij} = 1$, has a capacity of $u_{ij} \ge 0$ on the flow of all commodities, and $0$ otherwise.

The information structure of the problem is as follows.
We have to design the network, by determining the decisions $y$, while the value of the random vector $\omega$, and thus demand, is unknown.
After the realisation of $\omega$ becomes known, we can satisfy demand by sending a flow of commodities through the network.
Let $A^k = \{ (i, j) \in A \mid j \ne O(k) \text{ and } i \ne D(k) \}$ be the set of admissible arcs for each commodity $k \in K$, and let $x_{ij}^k \ge 0$ denote the flow variable representing the flow of commodity $k \in K$ through arc $(i, j) \in A^k$.
Here, $x$ depends on $\omega$, but we suppress that dependency in the notation.
The capacity of each arc $(i, j) \in A$ is shared between all the commodities for which $(i, j)$ is admissible. 
The goal is to design the network, that is, determine $y$ with minimal costs, such that the probability that demand can be satisfied for all commodities $k \in K$ by sending flow of amount $d^k(\omega)$ through the network is at least $1 - \alpha$ for some fixed parameter $\alpha \in [0, 1]$.

Throughout we assume that the support $\Omega$ of $\omega$ consists of a finite number $S$ of scenarios.
We label the realisations associated with these scenarios $\omega_s$, for $s \in [S]$, each occurring with probability $p_s > 0$.
This assumption is not very restrictive, since continuous distributions are typically approximated by finite discrete distributions using a sample average approximation~\citep{kleywegt_sample_2002, pagnoncelli_sample_2009}.

We say that \emph{$y \in \B^{|A|}$ is feasible for scenario $s \in [S]$} if there exists a solution to the feasibility subproblem $SP(y, \omega_s)$.
The feasibility subproblem consists of finding some vector of flows $x = (x)_{ij}^k \ge 0$ for each commodity $k \in K$ on the arcs $(i, j) \in A^k$, satisfying the following system of equations:
\begin{subequations}
    \label{eq:subproblem}
    \begin{alignat}{2}
            (SP(y, \omega)) \qquad 
                        \sum_{j \in N_i^-} x_{ji}^k &= \sum_{j \in N_i^+} x_{ij}^k, & \quad \forall k \in K,\forall i \in N \setminus \{ O(k), D(k) \}, & \qquad \text{(balance)}  \label{eq:balance} \\
        \sum_{k \in K \mid (i, j) \in A^k} x^k_{ij} &\le u_{ij}y_{ij},              & \quad \forall (i, j) \in A,                                       & \qquad \text{(capacity)} \label{eq:capacity} \\
                \sum_{j \in N_{D(k)}^-} x_{jD(k)}^k &\ge d^k(\omega),               & \quad \forall k \in K,                                            & \qquad \text{(demand)}   \label{eq:demand} \\
                                           x^k_{ij} &\ge 0,                         & \quad \forall k \in K, \forall (i, j) \in A^k,                    &                          \label{eq:box}    
    \end{alignat}
\end{subequations}
where $N^+_i = \{j \in N \mid (i, j) \in A\}$, and $N^-_i = \{j \in N \mid (j, i) \in A\}$.
The sets $A^k$ exclude arcs flowing into, or out of, each commodity's origin and destination nodes, respectively, which makes the formulation work by breaking any potential cycles around the origin and destination nodes.

Subproblem $SP(y, \omega)$ is allowed to be infeasible with a probability of at most $\alpha$ over the scenarios:
\begin{equation}
    \min_{y \in \B^{|A|}} \left\{ f^\top y~\middle|~\sum_{s \in [S]} p_s \1{SP(y, \omega_s) \text{ is infeasible}} \le \alpha \label{eq:math_problem} \right\}.
\end{equation}
Formulation~\eqref{eq:math_problem} is the type of model we aim to efficiently solve in this paper.
We assume throughout that the feasible region of~\eqref{eq:math_problem} is non-empty.
We next present a solution approach for this formulation in~\cref{sec:solution_approach}.

\section{Solution approach}
\label{sec:solution_approach}

We use the following, well-known result to reformulate~\eqref{eq:math_problem} as a tractable optimisation problem.
\begin{theorem}
    \label{th:sp_finite_cuts}
    Fix scenario $\bar s \in [S]$.
    The subset $Y_{\bar s} = \{ y \in \B^{|A|} \mid SP(y, \omega_{\bar s}) \text{ is feasible} \}$ can be described by finitely many linear inequalities (or \emph{feasibility cuts}) of the form:
    \begin{equation}
        0 \ge \gamma_{r\bar s} - \beta_{r\bar s}^\top y \qquad \forall r \in C_{\bar s},
        \label{eq:linear_ineq_cut_form}
    \end{equation}
    where $C_{\bar s}$ is the (index) set of feasibility cuts, $\gamma_{r\bar s}$ a scalar, and $\beta_{r\bar s}$ a vector in $\R^{|A|}$.
\end{theorem}
\proof{Proof.}
Fix $\bar y \in \B^{|A|}$, $\bar s \in [S]$, and assume $\bar y$ is not feasible for scenario $\bar s$.
Observe that $SP(\bar y, \omega_{\bar s})$ can be seen as a linear programme with a null minimisation objective.
This primal linear programme is feasible if and only if the objective of its dual, given by
\begin{subequations}
    \label{eq:dual_subproblem}
    \begin{alignat}{2}
        \max_{\substack{\mu~\text{free} \\ \pi, \lambda \ge 0}} \biggl\{ \sum_{k \in K} d^k(\omega) \lambda^k - \sum_{(i, j) \in A} u_{ij} y_{ij} \pi_{ij} \biggm|
            \mu_{O(k)}^k &= 0,                    & \quad \forall k \in K,                        & \label{eq:dual_orig} \\
            \mu_{D(k)}^k &= \lambda^k,            & \quad \forall k \in K,                        & \label{eq:dual_dest} \\
                \pi_{ij} &\ge \mu_j^k -  \mu_i^k, & \quad \forall k \in K,~\forall (i, j) \in A^k & \label{eq:dual_edge} \biggr\},
    \end{alignat}
\end{subequations}
is bounded by zero
Here, $\mu_i^k$, $\pi_{ij}$, and $\lambda^k$ are dual variables corresponding to the constraints in~\eqref{eq:balance},~\eqref{eq:capacity}, and~\eqref{eq:demand}, respectively, with $\mu_{O(k)}^k = 0$ and $\mu_{O(k)}^k = \lambda^k$ additionally defined for each commodity $k \in K$ for notational convenience.
Note that the feasible region of the dual programme~\eqref{eq:dual_subproblem} does not depend on a $y \in \B^{|A|}$.
In fact, this feasible region is a convex cone with finitely many extreme rays, and the dual objective is bounded by zero if the objective corresponding to all these extreme rays is bounded by zero.
Hence, every such extreme ray $(\mu^*, \pi^*, \lambda^*)$ defines a feasibility cut in~\eqref{eq:linear_ineq_cut_form} via
\[ \sum_{k \in K} d^k(\omega_{\bar s}) \lambda^{k*} - \sum_{(i, j) \in A} u_{ij} y_{ij} \pi_{ij}^* \le 0,\]
and together these inequalities fully describe the subset $Y_{\bar s} \subseteq \B^{|A|}$ that is feasible for $\bar s$.
\QED
\endproof

\begin{subequations}
    We use~\cref{th:sp_finite_cuts} to reformulate~\eqref{eq:math_problem} using a knapsack-like formulation~\citep{raike_dissection_1970,luedtke_branch-and-cut_2014}.
    In particular, we introduce binary variables $z_s \in \B$ for each scenario $s \in [S]$, deciding whether $y$ must be made feasible for scenario $s$ ($z_s = 0$) or not ($z_s = 1$).
    The knapsack constraint is given by
    \begin{equation}
        \sum_{s \in [S]} p_s z_s \le \alpha.
        \label{eq:knapsack_scenarios}
    \end{equation}
    Using the variables $z_s$ to determine whether $y$ should be made feasible for scenario $s$, we also impose the following feasibility cuts of the form~\eqref{eq:linear_ineq_cut_form}:
    \begin{equation}
        M_{rs} z_s \ge \gamma_{rs} - \beta_{rs}^\top y, \qquad \forall s \in [S],~r \in C_s,
        \label{eq:all_cuts}
    \end{equation}
    where each $M_{rs} > 0$ is sufficiently large to ensure the constraints are non-binding for $y$ when $z_s = 1$.
\end{subequations}

Together, constraints~\eqref{eq:knapsack_scenarios} and~\eqref{eq:all_cuts} are equivalent to the probabilistic constraint in formulation~\eqref{eq:math_problem}.
The difficulty in solving this reformulation lies in the number of linear inequalities that make up the constraints~\eqref{eq:all_cuts}: although finite due to~\cref{th:sp_finite_cuts} and $S < \infty$, in practice, the constraint set~\eqref{eq:all_cuts} is very large.
Instead of directly solving a model with all such constraints, we devise an iterative scheme where in each iteration we separate inequalities of the form of~\eqref{eq:all_cuts} from the infeasible scenarios.
Specifically, we propose the following master problem in iteration $l > 0$:
\begin{subequations}
    \label{eq:master_problem}
    \begin{alignat}{2}
        (MP_l) \quad \min_{\substack{y \in \B^{|A|} \\ z \in \B^S}} \biggl\{ f^\top y \biggm|
        \sum_{s \in [S]} p_s z_s &\le \alpha,                          &                 &                                         \label{eq:select_scenarios} \\
                      M_{rs} z_s &\ge \gamma_{rs} - \beta_{rs}^\top y, & \quad \forall r &\in [l - 1],~\forall s \in S_r \biggr\}. \label{eq:feas_cut}
    \end{alignat}  
\end{subequations}

Constraint~\eqref{eq:feas_cut} lists the feasibility cuts from earlier iterations $1 \le r < l$.
More precisely, the set $S_r \subseteq [S]$ represents the set of infeasible scenarios that were not allowed to be infeasible (w.r.t.\ the optimal solution $(y^*_r, z^*_r)$ to $MP_r$) in earlier iterations $r$ for which feasibility cuts were inserted into the master problem.

The rest of this section covers the details of our iterative scheme.
First, in~\cref{subsec:feasibility_cuts}, we discuss four alternative formulations of the subproblem $SP(y, \omega)$ from which to derive the values of $M$, $\gamma$, and $\beta$ that make up the feasibility cuts~\eqref{eq:feas_cut}.
Three of these formulations are commonly used in the literature, and we also present a novel formulation.
Second, when only a few iterations have passed, the feasibility cuts~\eqref{eq:feas_cut} derived for each scenario $s \in [S]$ are still a poor approximation of the feasible set.
As such, the master problem solutions are likely to be poor initially.
In~\cref{subsec:valid_inequalities} we present several valid inequalities that strengthen the master problem, particularly in these early iterations.
Such inequalities significantly speed-up convergence.

\subsection{Feasibility cuts}
\label{subsec:feasibility_cuts}

Given a solution $y \in \B^{|A|}$ and scenario $s \in S$, we want to efficiently determine if $y$ is feasible for scenario $s$.
An optimisation perspective is helpful in achieving this.
Observe that $SP(y, \omega_s)$ can be seen as a linear programme with a null minimisation objective.
If no flow vector $x$ satisfying constraints~\eqref{eq:subproblem} exists, the linear programme is infeasible.
Its dual programme given in~\eqref{eq:dual_subproblem} is then unbounded, and one can identify an extreme ray along which the dual objective can be arbitrarily improved.
A constraint based on this extreme ray can be added to the master problem as a Benders' feasibility cut, as was proposed originally by \cite{benders_partitioning_1962}.
Rather than identifying extreme rays, it is common practice to slightly alter the constraint system of the dual programme instead to avoid unboundedness altogether.
The feasibility cuts are then derived by guaranteeing that the dual programme's objective function is bounded by zero at every vertex of the dual feasible region (see \textit{e.g.}~\cite{minoux_networks_1989} for an early application of this idea to multicommodity network flow problems).

In this section we propose four such alternative formulations: first, our novel formulation for this specific dual programme that we term \textit{FlowMIS}, and then three general formulations from the literature that we apply to this dual programme.
The feasibility cuts derived from these alternative formulations have different properties, which we investigate in the single-commodity setting in~\cref{subsubsec:cut_properties_single_commodity}.

\paragraph{FlowMIS.}
The dual programme given by~\eqref{eq:dual_subproblem} is potentially unbounded because $\lambda^k$ can be arbitrarily increased, for each commodity $k \in K$.
Adding the restriction
\begin{equation}
    \sum_{k \in K} \lambda^k \le 1
    \label{eq:flowmis}
\end{equation}
to~\eqref{eq:dual_subproblem} is sufficient to avoid unboundedness.
This is equivalent to artificially lowering demand via a slack variable $t \ge 0$ in constraints~\eqref{eq:demand_flowmis_primal} in the primal formulation associated with the dual programme given by~\eqref{eq:dual_subproblem} and~\eqref{eq:flowmis}: 
\begin{subequations}
    \label{eq:flowmis_primal}
    \begin{alignat}{3}
        \min_{x, t \ge 0} && t                                           &                               &                                                                   & \notag \\
              \text{s.t.} &&                 \sum_{j \in N_i^-} x_{ji}^k &= \sum_{j \in N_i^+} x_{ij}^k, & \quad \forall k \in K,\forall i \in N \setminus \{ O(k), D(k) \}, & \label{eq:balance_flowmis_primal} \\
                          && \sum_{k \in K \mid (i, j) \in A^k} x^k_{ij} &\le u_{ij}y_{ij},              & \quad \forall (i, j) \in A,                                       & \label{eq:capacity_flowmis_primal} \\
                          &&         \sum_{j \in N_{D(k)}^-} x_{jD(k)}^k &\ge d^k(\omega) - t,           & \quad \forall k \in K.                                            & \label{eq:demand_flowmis_primal}
    \end{alignat}
\end{subequations}
We term this formulation \textit{FlowMIS}, since it works for primal problems with a network flow structure by modifying only the demand constraints. 

\paragraph{SNC.}
This formulation adds a constraint that forces the sum of all non-free dual variables to lie in the interval $[0, 1]$.
That constraint was first discussed by~\cite{balas_mixed_1996} and \cite{balas_modified_1997}, and~\cite{fischetti_separation_2011} termed it the \emph{standard normalisation condition} (SNC).
In particular, SNC adds the following constraint to~\eqref{eq:dual_subproblem}:
\begin{equation}
    \sum_{(i, j) \in A} \pi_{ij} + \sum_{k \in K} \lambda^k \le 1.
    \label{eq:snc}
\end{equation}
This dual constraint corresponds to a single slack variable inserted in each of the capacity and demand constraints in the primal formulation, and penalising its value in the objective.
Compared to the FlowMIS primal of~\eqref{eq:flowmis_primal}, the SNC formulation also inserts the slack variable $t$ into the capacity constraints, replacing~\eqref{eq:capacity_flowmis_primal} with
\[ \sum_{k \in K \mid (i, j) \in A^k} x^k_{ij} \le u_{ij}y_{ij} + t, \quad \forall (i, j) \in A. \]

\paragraph{MIS.}
\cite{fischetti_note_2010} look directly at the implication of~\eqref{eq:snc} on the primal constraint system.
Finding constraints that require slack variables is a search for a \emph{minimal infeasible subsystem} (MIS) of constraints that together make $SP(y, \omega)$ infeasible.
They observe that constraints that do not involve $y$ are static conditions that are always active.
Since those constraints are never part of an infeasible subset of constraints, there is no reason to insert a slack variable in them.
As such, they propose the following refinement of~\eqref{eq:snc}:
\begin{equation}
    \sum_{(i, j) \in A} \pi_{ij} \le 1,
    \label{eq:mis}
\end{equation}
since only capacity constraints involving arcs in $A$ depend on $y$.
Compared to the FlowMIS primal of~\eqref{eq:flowmis_primal}, the MIS formulation only inserts the slack variable $t$ into the capacity constraints, replacing~\eqref{eq:capacity_flowmis_primal} and~\eqref{eq:demand_flowmis_primal} with
\begin{alignat*}{2}
    \sum_{k \in K \mid (i, j) \in A^k} x^k_{ij} &\le u_{ij}y_{ij} + t,              & \quad \forall (i, j) \in A, & \\
            \sum_{j \in N_{D(k)}^-} x_{jD(k)}^k &\ge d^k(\omega),               & \quad \forall k \in K.      &
\end{alignat*}
CPLEX uses the MIS formulation internally for its automatic Benders' decomposition feature~\citep{bonami_implementing_2020}.

\paragraph{BB.}
Finally, the \emph{basic Benders} (BB) formulation ensures each non-free variable is bounded to the unit interval.
In particular, BB adds the following constraints to~\eqref{eq:dual_subproblem}:
\begin{equation}
    \label{eq:bb}
    \begin{aligned}
         \pi_{ij} &\le 1, \qquad \forall (i, j) \in A, \\
        \lambda^k &\le 1, \qquad \forall k \in K.
    \end{aligned}
\end{equation}
These dual constraints correspond to inserting unique slack variables in each of the capacity~\eqref{eq:capacity} and demand~\eqref{eq:demand} constraints in the primal, and penalising their values in the objective.
We term this formulation `basic Benders' because relaxing each inequality with a unique slack variable is a straightforward way of dealing with primal infeasibility~\citep{birge_two-stage_2011}.
Compared to the FlowMIS primal of~\eqref{eq:flowmis_primal}, the BB has many more slack variables (one for each capacity and demand constraint).

\paragraph{Determining $M$, $\beta$, and $\gamma$.}
Let $y^*_r \in \B^{|A|}$ be the solution to $MP_r$~\eqref{eq:master_problem} in iteration $r > 0$.
Assume $y^*_r$ is not feasible for a fixed scenario $\bar s \in [S]$.
After solving any of the four alternative formulations of dual programme~\eqref{eq:dual_subproblem} with $y^*_r$ and $\omega_{\bar s}$, a feasibility cut for $MP_r$ is derived as follows.
Let $\mu^*$, $\pi^*$, and $\lambda^*$  be the optimal dual solutions.
The optimal dual objective is then given by
\[ 0 \le \sum_{k \in K} d^k(\omega_{\bar s}) \lambda^{k*} - \sum_{(i, j) \in A} u_{ij} y^*_{ij} \pi_{ij}^* < \infty. \]
This objective value is bounded by construction, and strictly positive only when the original primal problem (without slack variables) is infeasible.
We then impose the following constraint on $y$:
\begin{equation}
            M_{r\bar s} z_{\bar s} \ge 
                \underbrace{\sum_{k \in K} d^k(\omega_{\bar s}) \lambda^{k*}}_{=\gamma_{r\bar s}} 
                - \underbrace{\sum_{(i, j) \in A} u_{ij} y_{ij} \pi_{ij}^*}_{=\beta_{r\bar s}^\top y}. \label{eq:dual_feas_cut}
\end{equation}
Constraint~\eqref{eq:dual_feas_cut} is clearly a feasibility cut of the form~\eqref{eq:feas_cut}.
Proposition~\ref{prop:M_is_gamma} shows we can select $\gamma_{r\bar s}$ as a tight value for $M_{r\bar s}$.
\begin{proposition}
    \label{prop:M_is_gamma}
    In a feasibility cut of the form~\eqref{eq:dual_feas_cut}, $M_{r\bar s} = \gamma_{r\bar s}$ is a tight value for $M_{r\bar s}$.
\end{proposition}
\proof{Proof.} 
    When $z_{\bar s} = 1$, we must have $M_{r\bar s} \ge \gamma_{r\bar s} - \beta_{r\bar s}^\top y$.
    We know that $u_{ij} \ge 0$ for all $(i, j) \in A$.
    An optimal solution $(\mu^*, \pi^*, \lambda^*)$ to any of the four (re)formulations of~\eqref{eq:dual_subproblem} satisfies $\pi^* \ge 0$.
    Since $y \ge 0$ as well, we have that $\beta_{r\bar s}^\top y \ge 0$.
    Thus $M_{r\bar s} \ge \gamma_{r\bar s} - \beta_{r\bar s}^\top y$ is always satisfied by choosing $M_{r\bar s} = \gamma_{r\bar s}$, and no smaller value would suffice.
\QED
\endproof

\paragraph{Metric inequalities.}
Feasibility cuts of the form~\eqref{eq:dual_feas_cut} can be strengthened further by raising the value of $\gamma_{r\bar s}$.
In particular, by solving $|K|$ shortest path problems, it is possible to derive a so-called \textit{metric inequality} that improves the value of $\gamma_{r\bar s}$ and strengthens the resulting feasibility cut.
To achieve this we implement the procedure of~\cite{costa_benders_2009}: for each commodity $k \in K$, we solve a shortest path problem between $O(k)$ and $D(k)$, where the arc weights are given by $\pi^*_{ij}$ for arcs $(i, j) \in A$.
This results in a shortest path $P^k \subseteq A^k$ for each commodity $k \in K$.
We now determine $\gamma_{r\bar s}$ as
\[ \gamma_{r\bar s} = \sum_{k \in K} d^k(\omega_{\bar s}) \sum_{(i, j) \in P^k} \pi^*_{ij}. \]

\subsubsection{Cut properties in the single-commodity case}
\label{subsubsec:cut_properties_single_commodity}

In this section, we discuss the differences in the feasibility cuts derived from the subproblem formulations that we have presented.
We do this for the single-commodity case, that is, for $|K| = 1$, since we can make use of the max-flow min-cut theorem in this case.
We will prove that each of the different feasibility subproblem formulations results in feasibility cuts that are different kinds of minimum-cuts in the network with capacities $u_{ij} \bar y_{ij}$ for $(i, j) \in A$.
In particular,~\cref{th:flowmis_single_commodity} shows that FlowMIS finds a minimum-capacity cut, whereas~\cref{th:snc_single_commodity,th:mis_single_commodity} show that the SNC and MIS formulations find cuts that are cardinality-constrained, respectively.
Such cardinality-constrained cuts are not, in general, the same as the minimum-capacity cuts that FlowMIS obtains.
Finally, although~\cref{cor:flowmis_opt_for_bb} shows that the BB formulation can also find minimum-capacity cuts, we will demonstrate in~\cref{sec:numerical_experiments} that its increased size results in substantially worse performance than any of the other three formulations.

\begin{theorem}
    \label{th:flowmis_single_commodity} 
    Assume $|K| = 1$.
    Fix $\bar y \in \B^{|A|}$ and $\bar s \in [S]$, and assume $\bar y$ is not feasible for $\bar s$.
    Let $C^*_{\text{FlowMIS}}$ be a $(O(1), D(1))$ cut in $G$ of minimum capacity.
    A feasibility cut derived from the FlowMIS formulation of~\eqref{eq:dual_subproblem} and~\eqref{eq:flowmis} is equivalent to
    \[ \sum_{(i, j) \in C^*_{\text{FlowMIS}}} u_{ij} y_{ij} \ge d^1 (\omega_{\bar s}). \]
\end{theorem}
\proof{Proof.}
Consider both the primal FlowMIS formulation in~\eqref{eq:flowmis_primal} and its dual in~\eqref{eq:dual_subproblem} and~\eqref{eq:flowmis}.
We will construct feasible primal and dual solutions to~\eqref{eq:flowmis_primal} and~\eqref{eq:dual_subproblem}--\eqref{eq:flowmis}, respectively, and prove that they are both optimal by showing that their objective values are the same.

First consider the primal FlowMIS formulation of~\eqref{eq:flowmis_primal}.
Here, the slack variable $t \ge 0$ decreases the demand $d^1(\omega_{\bar s})$.
A natural condition for the feasibility of~\eqref{eq:flowmis_primal} is that this $t$ should be selected such that all $(O(1), D(1))$ cuts in the graph $G$ have capacity of at least $d^1(\omega_{\bar s}) - t$.
Let $\mathcal{C}$ be the set of all $(O(1), D(1))$ cuts in $G$.
We thus must have that
\[ \sum_{(i, j) \in C} u_{ij} \bar y_{ij} \ge d^1(\omega_{\bar s}) - t \qquad \forall C \in \mathcal{C}. \]
The smallest $t$ that meets this condition is given by
\begin{equation}
    t^* \coloneqq \max_{C \in \mathcal{C}} \left\{ d^1(\omega_{\bar s}) - \sum_{(i, j) \in C} u_{ij} \bar y_{ij} \right\}.
    \label{eq:necessary_condition_flowmis}
\end{equation}
Since $\bar y$ is infeasible for $\bar s$, we must have that $t^* > 0$.
Let $C^*_{\text{FlowMIS}}$ be a cut that maximises $t$ in~\eqref{eq:necessary_condition_flowmis}.
Such a cut $C^*_{\text{FlowMIS}}$ has the lowest capacity of all $(O(1), D(1))$ cuts in $G$.
An optimal solution to the FlowMIS formulation of~\eqref{eq:dual_subproblem} and~\eqref{eq:flowmis} is constructed as follows.
Set $\lambda^{1*} = 1$, and take for $i \in N$
\[
    \mu_i^1 = \begin{cases}
        1 & \text{if there exists an $(i, D(1))$-path in } A \setminus C^*_{\text{FlowMIS}}, \\
        0 & \text{otherwise},
    \end{cases}
\]
and for $(i, j) \in A$,
\[
    \pi^*_{ij} = \begin{cases}
        1 & \text{if } (i, j) \in C^*_{\text{FlowMIS}}, \\
        0 & \text{otherwise}.
    \end{cases}
\]
This solution is clearly feasible, and attains the same objective value as the optimal primal solution:
\[
    \lambda^{1*} d^1(\omega_{\bar s}) - \sum_{(i, j) \in A} \pi^*_{ij} u_{ij} \bar y_{ij}
        = d^1(\omega_{\bar s}) - \sum_{(i, j) \in C^*_{\text{FlowMIS}}} \pi^*_{ij} u_{ij} \bar y_{ij} = t^*.
\]
The resulting feasibility cut is then given by
\[
    \lambda^{1*} d^1(\omega_{\bar s}) - \sum_{(i, j) \in A} \pi_{ij}^* u_{ij} y_{ij} \le 0 
        \implies \sum_{(i, j) \in C^*_{\text{FlowMIS}}} u_{ij} y_{ij} \ge d^1(\omega_{\bar s}),
\]
which concludes the proof.
\QED
\endproof

The FlowMIS formulation thus finds a minimum-capacity cut, and requires the capacity of this cut to be increased to at least $d^1 (\omega_{\bar s})$ when $y$ must be feasible for scenario $\bar s$.
The following theorems show that the SNC and MIS formulations may find cuts with a lower cardinality than the minimum-capacity cut.

\begin{theorem}
    \label{th:snc_single_commodity}
    Assume $|K| = 1$.
    Fix $\bar y \in \B^{|A|}$ and $\bar s \in [S]$, and assume $\bar y$ is not feasible for $\bar s$.
    Let $\mathcal{C}$ be the set of all $(O(1), D(1))$ cuts in $G$.
    Take
    \begin{equation}
        \label{eq:snc_cut_condition}
        C^*_{\text{SNC}} = \arg \max_{C \in \mathcal{C}} \left\{ \frac{d^1(\omega_{\bar s}) - \sum_{(i, j) \in C} u_{ij} \bar y_{ij}}{|C| + 1} \right\}.
    \end{equation}
    A feasibility cut derived from the SNC formulation of~\eqref{eq:dual_subproblem} and~\eqref{eq:snc} is equivalent to
    \[ \sum_{(i, j) \in C^*_{\text{SNC}}} u_{ij} y_{ij} \ge d^1(\omega_{\bar s}). \]
\end{theorem}
\proof{Proof.}
We proceed in much the same way as we did proving~\cref{th:flowmis_single_commodity}, by constructing optimal primal and dual solutions, and showing their objective values coincide.
Now, the primal of the dual SNC formulation in~\eqref{eq:dual_subproblem} and~\eqref{eq:snc} is given by
\begin{alignat*}{3}
    \min_{x, t \ge 0} && t                                   &                               &                                           & \\
          \text{s.t.} &&         \sum_{j \in N^-_i} x^1_{ji} &= \sum_{j \in N^+_i} x^1_{ij}, & \forall i \in N \setminus \{O(1), D(1)\}, & \\
                      &&                            x^1_{ij} &\le u_{ij} y_{ij} + t,         & \forall (i, j) \in A,                     & \\ 
                      && \sum_{j \in N^-_{D(k)}} x^1_{jD(k)} &\ge d^1(\omega_{\bar s}) - t.  &                                           &
\end{alignat*}
Here, the slack variable $t \ge 0$ increases the capacity of all edges, and decreases the demand.
A natural condition on feasibility is that this $t$ should be selected such that all $(O(1), D(1))$ cuts in the graph $G$ have capacity of at least $d^1(\omega_{\bar s}) - t$.
Let $\mathcal{C}$ be the set of all $(O(1), D(1))$ cuts in $G$.
We thus must have that
\[ \sum_{(i, j) \in C} u_{ij} \bar y_{ij} + t|C| \ge d^1(\omega_{\bar s}) - t \qquad \forall C \in \mathcal{C}. \]
The smallest $t$ that meets this condition is given by
\begin{equation}
    t^* \coloneqq \max_{C \in \mathcal{C}} \left\{ \frac{d^1(\omega_{\bar s}) - \sum_{(i, j) \in C} u_{ij} \bar y_{ij}}{|C| + 1} \right\}.
    \label{eq:necessary_condition_snc}
\end{equation}
Since $\bar y$ is infeasible for $\bar s$, we must have that $t^* > 0$.
Let $C^*_{\text{SNC}}$ be a cut that maximises $t$ in~\eqref{eq:necessary_condition_snc}.
An optimal solution to the SNC formulation of~\eqref{eq:dual_subproblem} and~\eqref{eq:snc} is constructed as follows.
Set $\lambda^{1*} = \frac{1}{|C^*| + 1}$, and take for $i \in N$
\[
    \mu_i^1 = \begin{cases}
        \frac{1}{|C^*| + 1} & \text{if there exists an $(i, D(1))$-path in } A \setminus C^*_{\text{SNC}}, \\
                          0 & \text{otherwise},
    \end{cases}
\]
and for $(i, j) \in A$,
\[
    \pi^*_{ij} = \begin{cases}
        \frac{1}{|C^*| + 1} & \text{if } (i, j) \in C^*_{\text{SNC}}, \\
                          0 & \text{otherwise}.
    \end{cases}
\]
This dual solution is clearly feasible, and attains the same objective value as the primal solution:
\begin{align*}
    \lambda^{1*} d^1(\omega_{\bar s}) - \sum_{(i, j) \in A} \pi_{ij}^* u_{ij} \bar y_{ij} 
        &= \frac{1}{|C^*_{\text{SNC}}| + 1} d^1(\omega_{\bar s}) - \frac{1}{|C^*_{\text{SNC}}| + 1} \sum_{(i, j) \in C^*_{\text{SNC}}} u_{ij} \bar y_{ij} \\
        &= \max_{C \in \mathcal{C}} \left\{ \frac{d^1 (\omega_{\bar s}) - \sum_{(i, j) \in C} u_{ij} \bar y_{ij}}{|C| + 1} \right\} = t^*.
\end{align*}
The resulting feasibility cut is then given by
\[
    \lambda^{1*} d^1(\omega_{\bar s}) - \sum_{(i, j) \in A} \pi_{ij}^* u_{ij} y_{ij} \le 0
        \implies \sum_{(i, j) \in C^*_{\text{SNC}}} u_{ij} y_{ij} \ge d^1(\omega_{\bar s}),
\]
which concludes the proof.
\QED
\endproof

The SNC formulation may thus find cuts that do not have minimum capacity, but rather have a large capacity shortfall over cardinality ratio.
Such cuts may also be minimum capacity cuts, but this is not generally the case due to the cut cardinality aspect.
Also for SNC, the derived feasibility constraint requires the capacity of this cut to be increased to at least $d^1(\omega_{\bar s})$ when $y$ must be feasible for scenario $\bar s$.

\begin{theorem}
    \label{th:mis_single_commodity}
    Assume $|K| = 1$.
    Fix $\bar y \in \B^{|A|}$ and $\bar s \in [S]$, and assume $\bar y$ is not feasible for $\bar s$.
    Let $\mathcal{C}$ be the set of all $(O(1), D(1))$ cuts in $G$.
    Take
    \begin{equation}
        \label{eq:mis_cut_condition}
        C^*_{\text{MIS}} = \arg \max_{C \in \mathcal{C}} \left\{ \frac{d^1(\omega_{\bar s}) - \sum_{(i, j) \in C} u_{ij} \bar y_{ij}}{|C|} \right\}.
    \end{equation}
    A feasibility cut derived from the MIS formulation of~\eqref{eq:dual_subproblem} and~\eqref{eq:mis} is equivalent to
    \[ \sum_{(i, j) \in C^*_{\text{MIS}}} u_{ij} y_{ij} \ge d^1(\omega_{\bar s}). \]
\end{theorem}
\proof{Proof.} 
Analogous to the proof of \cref{th:snc_single_commodity}.
\QED
\endproof

Like the SNC formulation, the MIS formulation finds a type of cardinality-constrained cut.
It similarly requires that the capacity of this cut should be increased to at least $d^1(\omega_{\bar s})$ when $y$ must be feasible for scenario $\bar s$.
The main difference between the SNC and MIS formulations is in the denominator of the cut condition: since the SNC constrains both the $\pi$ and $\lambda$ dual variables in~\eqref{eq:snc}, the denominator in~\eqref{eq:snc_cut_condition} is the cut cardinality plus one, unlike the MIS formulation, which does not constrain the $\lambda$ dual variables in~\eqref{eq:mis} and thus divides by the cut cardinality in~\eqref{eq:mis_cut_condition}.

The following corollary, which states that an optimal solution to the FlowMIS formulation is also optimal for BB.
Thus, the BB formulation is also able to find minimum-capacity cuts.
\begin{corollary}
    \label{cor:flowmis_opt_for_bb}
    Let $|K| = 1$.
    An optimal dual solution $(\mu^*, \pi^*, \lambda^*)$ to the FlowMIS formulation of~\eqref{eq:dual_subproblem} and~\eqref{eq:flowmis} is also optimal for the BB formulation of~\eqref{eq:dual_subproblem} and~\eqref{eq:bb}.
\end{corollary}
\proof{Proof.}
For $|K| = 1$, the dual feasible region corresponding to BB is a subset of that corresponding to FlowMIS, since FlowMIS only imposes the additional constraint $\lambda^1 \le 1$, whereas BB in addition also imposes the constraints $\pi_{ij} \le 1$ for all $(i, j) \in A$.
However, notice that each optimal FlowMIS solution $(\mu^*, \pi^*, \lambda^*)$ as defined in the proof of~\cref{th:flowmis_single_commodity} satisfies these additional constraints.
Hence, this solution is also optimal for the BB formulation of ~\eqref{eq:dual_subproblem} and~\eqref{eq:bb}.
\QED
\endproof

\Cref{cor:flowmis_opt_for_bb} suggests BB might perform about the same as FlowMIS, but that is unlikely to be the case in practice.
First, the argument requires $|K| = 1$, since only then the feasible regions of BB and FlowMIS are the same.
Second, the feasible region of BB contains many more vertices than that of FlowMIS, by virtue of the additional constraints on $\pi_{ij}$.
Although the FlowMIS solution is an optimal vertex, there are very likely to be other vertices that are also optimal, which result in different feasibility cuts.
We will compare the numerical performance of these two formulations for single-commodity instances in~\cref{subsec:results_singlecommodity}.

\subsection{Valid inequalities}
\label{subsec:valid_inequalities}

The master problem $MP$ of~\eqref{eq:master_problem} contains little information about the structure of the network design problem, since all problem-specific constraints are present in the subproblem $SP$ of~\eqref{eq:subproblem}.
Its solutions will thus be poor for early iterations when few feasibility cuts are present, which results in slow convergence.
In this subsection we develop a scenario creation strategy in the manner of~\cite{crainic_partial_2021} to strengthen the master problem $MP$.
The additional structure this scenario provides helps the master problem find better solutions earlier, so that the decomposition algorithm needs fewer iterations to converge.

We first detail our scenario creation strategy, and then present~\cref{prop:marginal_scenario_is_valid} to establish its validity for $MP$.
For each commodity $k \in K$, let
\[ \bar d^k = \min_{z \in \B^S} \left\{ \sum_{s = 1}^S d^k(\omega_s) p_s (1 - z_s) \biggm| \sum_{s = 1}^S p_s z_s \le \alpha \right\}. \]
In general, the minimisation problem to determine $\bar d^k$ can be interpreted as a knapsack problem in which the scenarios $[S]$ correspond to the items, the probabilities $p_s$ to their weights, and $\alpha$ to the size of the knapsack.
However, if all probabilities are equal, that is, if $p_s = \frac{1}{S}$ for all $s \in [S]$, then we may obtain $\bar d^k$ for every $k \in K$ by letting $\sigma_1^k, \ldots, \sigma_S^k$ denote a permutation of the scenarios such that $d^k(\omega_{\sigma_1^k}) \le d^k(\omega_{\sigma_2^k}) \le \ldots \le d^k(\omega_{\sigma_S^k})$, and defining
\[ \bar d^k = \frac{1}{S} \sum_{s = 1}^{\lceil (1 - \alpha) S \rceil} d^k(\omega_{\sigma_s^k}). \]
Then, we propose adding the following constraints and variables $\bar x$ to $MP$ to strengthen the formulation:
\begin{subequations}
    \label{eq:vi_scenario}
    \begin{alignat}{2}
                        \sum_{j \in N_i^-} \bar x_{ji}^k &= \sum_{j \in N_i^+} \bar x_{ij}^k, & \quad \forall k \in K,\forall i \in N \setminus \{ O(k), D(k) \}, & \label{eq:vi_scen_balance} \\
        \sum_{k \in K \mid (i, j) \in A^k} \bar x^k_{ij} &\le u_{ij}y_{ij},                   & \quad \forall (i, j) \in A,                                       & \label{eq:vi_scen_capacity} \\
                \sum_{j \in N_{D(k)}^-} \bar x_{jD(k)}^k &\ge \bar d^k,                       & \quad \forall k \in K,                                            & \label{eq:vi_scen_demand} \\
                                           \bar x^k_{ij} &\ge 0,                              & \quad \forall k \in K, \forall (i, j) \in A^k.                    & \label{eq:vi_scen_box}
    \end{alignat}
\end{subequations}

\begin{proposition}
    \label{prop:marginal_scenario_is_valid}
    Constraints~\eqref{eq:vi_scen_balance}--\eqref{eq:vi_scen_box} are valid for $MP$.
\end{proposition}
\proof{Proof.}
Consider the deterministic equivalent formulation of $MP$ and all subproblems $SP(y, \omega_s)$ for $s \in [S]$, given by
\begin{subequations}
    \label{eq:deq}
    \begin{alignat}{3}
    \min_{x, y, z} && \sum_{(i, j) \in A} f_{ij} y_{ij},             &                                  &                                                                                      & \label{eq:deq_objective} \\
       \text{s.t.} &&                 \sum_{j \in N_i^-} x_{ji}^{ks} &= \sum_{j \in N_i^+} x_{ij}^{ks}, & \quad \forall k \in K,\forall i \in N \setminus \{ O(k), D(k) \},\forall s \in [S],  & \label{eq:deq_scen_balance} \\
                   && \sum_{k \in K \mid (i, j) \in A^k} x^{ks}_{ij} &\le u_{ij}y_{ij},                 & \quad \forall (i, j) \in A,\forall s \in [S],                                        & \label{eq:deq_scen_capacity} \\
                   &&         \sum_{j \in N_{D(k)}^-} x_{jD(k)}^{ks} &\ge d^{k}(\omega_s) (1 - z_s),    & \quad \forall k \in K,\forall s \in [S],                                             & \label{eq:deq_scen_demand} \\
                   &&                       \sum_{s \in [S]} p_s z_s &\le \alpha,                       & \quad \forall k \in K,\forall s \in [S],                                             & \label{eq:deq_scen_limit} \\
                   &&                                    x^{ks}_{ij} &\ge 0,                            & \quad \forall k \in K, \forall (i, j) \in A^k,\forall s \in [S],                     & \label{eq:deq_scen_box1} \\
                   &&                                         y_{ij} &\in \B,                           & \quad \forall (i, j) \in A,                                                          & \label{eq:deq_scen_box2} \\
                   &&                                            z_s &\in \B,                           & \quad \forall s \in [S].                                                             & \label{eq:deq_scen_box3}
    \end{alignat}
\end{subequations}
From the constraints in~\eqref{eq:deq_scen_balance},~\eqref{eq:deq_scen_capacity},~\eqref{eq:deq_scen_demand} and~\eqref{eq:deq_scen_limit}, we obtain valid inequalities by taking a weighted sum over the scenarios, where the weights correspond to the probabilities $p_s$ for each scenario $s \in [S]$.
The resulting valid inequalities are
\begin{alignat}{2}
                    \sum_{j \in N_i^-} \bar x_{ji}^{ks} &= \sum_{j \in N_i^+} \bar x_{ij}^{ks}, & \quad \forall k \in K,\forall i \in N \setminus \{ O(k), D(k) \},\forall s \in [S],  & \notag \\
    \sum_{k \in K \mid (i, j) \in A^k} \bar x^{ks}_{ij} &\le u_{ij}y_{ij},                      & \quad \forall (i, j) \in A,\forall s \in [S],                                        & \notag \\
            \sum_{j \in N_{D(k)}^-} \bar x_{jD(k)}^{ks} &\ge d^{k}(\omega_s) (1 - z_s),         & \quad \forall k \in K,\forall s \in [S],                                             & \label{eq:vi_proof_rhs} \\
                               \sum_{s \in [S]} p_s z_s &\le \alpha,                            & \quad \forall k \in K,\forall s \in [S],                                             & \notag \\
                                       \bar x^{ks}_{ij} &\ge 0,                                 & \quad \forall k \in K, \forall (i, j) \in A^k,\forall s \in [S],                     & \notag
\end{alignat}
where $\bar x_{ij}^k = \sum_{s \in [S]} p_s x_{ij}^{ks}$ for all $k \in K$ and $(i, j) \in A^k$.
The desired result follows by observing that for every $z \in \B^S$ satisfying~\eqref{eq:deq_scen_demand} and for every commodity $k \in K$, $\bar d^k$ is a lower bound for the right-hand side in~\eqref{eq:vi_proof_rhs}.
\QED
\endproof


\vspace{1em}

Our treatment of relevant valid inequalities concludes the presentation of our solution approach to efficiently solve problems of the form~\eqref{eq:math_problem}.
Next, in~\cref{sec:numerical_experiments}, we show the effectiveness of our approach on a set of numerical experiments.

\section{Numerical experiments}
\label{sec:numerical_experiments}

In this section we numerically illustrate the performance of our solution approach of \cref{sec:solution_approach}.
We first briefly discuss some implementation details in~\cref{subsec:implementation_details}, present a brief experimental design to generate benchmark instances in~\cref{subsec:exp_design}, and finally present an analysis of the benchmark solutions in~\cref{subsec:results_singlecommodity} and~\cref{subsec:results_multicommodity}.

Our solution approach is implemented in Python 3.10 using Gurobi 10.0.
\makeatletter
\if@BLINDREV
The implementation is freely available at \texttt{URL redacted for double blind review}, under a liberal MIT license.
\else
The implementation is freely available at \url{https://github.com/N-Wouda/CC-NDP}, under a liberal MIT license.
\fi
\makeatother
Here one may also find the experimental instances generated in~\cref{subsec:exp_design}.
Unless otherwise noted, each experimental instance is solved using a single core of a 2.45GHz AMD 7763 processor, 32GB of memory, and two hours of run-time.
We fix $\alpha = 0.1$ throughout this section.

\subsection{Implementation details}
\label{subsec:implementation_details}

We apply our iterative scheme inside the branch-and-bound tree of the master problem, using the callback functionality that is available in most modern solvers.
In particular, the callback function is called every time Gurobi finds a new incumbent solution $(\hat y, \hat z)$.
For every $s \in [S]$, if $\hat z_s = 0$, we solve one of the subproblem formulations of~\cref{subsec:feasibility_cuts} to determine if $\hat y$ is feasible for $s$.
If it is not, we derive a feasibility cut from the subproblem solution, and add that to the master problem as a so-called lazy cuts.
Our implementation thus adds many cuts simultaneously, potentially one for each of the $S$ scenarios.
This does not slow down the master problem in practice, as Gurobi's cut pool manager applies only those lazy cuts that are relevant.

For exposition and consistency with the existing literature, we presented the subproblems from a dual perspective in~\cref{subsec:feasibility_cuts}.
In our implementation we implement the primal formulation, that is, a suitable modification of $SP$ directly, not of its dual.
Since only the right-hand sides of $SP(y, \omega)$ change for different vectors $y$ and fixed $\omega$, we solve all subproblem formulations using the dual simplex algorithm.
This lets us benefit from warm-starting at the previous solution for an earlier vector $y'$, if such a solution exists.

\subsection{Experimental design}
\label{subsec:exp_design}

We use a subset of the well-known \textbf{R} instances as a basis for our benchmarks.
These instances are widely used in the literature, see for example~\cite{boland_proximity_2016},~\cite{chouman_commodity_2017},~\cite{rahmaniani_accelerating_2018}, and~\cite{crainic_exact_2021}.
The \textbf{R} instances follow the classification scheme \textbf{R}$x$-$y$.
We consider the instance groups with `$x$' in $\{04, 05, \ldots, 10\}$, and focus on instances with tight capacity constraints (those with `$y$', the fixed cost/capacity parameter, in $\{7, 8, 9\}$).
These instances provide the network topology, arc capacities and fixed costs, and commodities.
The commodity demands are taken from the corresponding $1,000$ (uncorrelated) demand scenarios of~\cite{rahmaniani_asynchronous_2024} for each \textbf{R} instance.
We generate sets of $16$, $32$, $64$, $128$, $256$, and $512$ scenarios from these $1,000$ scenarios: the first 16 go into the first instance, the next 32 into the second, and so on until the $512$ scenario instance is generated.
This results in a total of $126$ instances.
\Cref{tab:instances} summarises the number of nodes, arcs, commodities, and scenarios in each of the instances.
We label each instance R$x$-$y$-$z$, where $x \in \{04, 05, \ldots, 10\}$, $y \in \{7, 8, 9\}$, and $z \in \{16, 32, 64, 128, 256, 512\}$.

\begin{table}
    \caption{Number of nodes, arcs, commodities, and scenarios in each instance group.}
    \label{tab:instances}

    \centering

    \begin{tabular}{lrrrl}
    \toprule
     & $|N|$ & $|A|$ & $|K|$ & $S$ \\
    \midrule
    R04 & 10 & 60 & 10 & \{16, 32, 64, 128, 256, 512\} \\
    R05 & 10 & 60 & 25 & \{16, 32, 64, 128, 256, 512\} \\
    R06 & 10 & 60 & 50 & \{16, 32, 64, 128, 256, 512\} \\
    R07 & 10 & 82 & 10 & \{16, 32, 64, 128, 256, 512\} \\
    R08 & 10 & 83 & 25 & \{16, 32, 64, 128, 256, 512\} \\
    R09 & 10 & 83 & 50 & \{16, 32, 64, 128, 256, 512\} \\
    R10 & 20 & 120 & 40 & \{16, 32, 64, 128, 256, 512\} \\
    \bottomrule
    \end{tabular}
\end{table}

We also construct 126 \textit{single-commodity} instances to explore the implications of the properties we proved in~\cref{subsubsec:cut_properties_single_commodity}.
These instances all have the same attributes as the 126 multi-commodity instances whose generation we just explained, except that only the first commodity is retained: all other commodities are removed from the instance.

\subsection{Results on single-commodity instances}
\label{subsec:results_singlecommodity}

\Cref{tab:single-commodity} presents the results of solving the single-commodity instances, with the valid inequalities of~\cref{subsec:valid_inequalities} enabled.
We do not use the metric inequality strengthening for the single-commodity instances since it is not effective when $|K| = 1$. 
We observe that FlowMIS on average outperforms the other formulations, although some variation exists between the various instance groups (particularly R10, which MIS and SNC both solve significantly faster).
The FlowMIS, MIS, and SNC formulation solve all instances within the two hour time limit, but BB did not manage to solve a single large R09 instance in time.
FlowMIS solves 48 instances the fastest, followed by SNC with 35, MIS with 34, and finally BB with just 9.
The FlowMIS, MIS, and SNC formulations are all more than eight times faster than BB on average.

BB performs considerably worse than FlowMIS, even though an optimal FlowMIS solution is also optimal for BB by virtue of~\cref{cor:flowmis_opt_for_bb}.
It is clear from~\cref{tab:single-commodity} that they do not find the same feasibility cuts, as the BB formulation requires two to three times more iterations to solve instances than the FlowMIS formulation needs.
Investigating the feasibility cuts returned by both formulations, we find that the cuts created from the BB formulation are much less sparse: the $\beta$ vector contains many more non-zero elements than a cut created from the FlowMIS formulation does.
Upon closer inspection, the reason appears to be that BB sets more $\pi_{ij}$ variables to $1$ when those $\pi_{ij}$ variables are associated with arcs $(i, j) \in A$ for which the construction decision $y_{ij} = 0$.
This choice does not affect optimality, but it does result in cuts with increased cardinality, and that affects the numerical performance.

\begin{table}
    \caption{
        Average number of iterations, run-time, and number of solved single-commodity instances (`\#S') for each subproblem formulation, by instance group, with valid inequalities enabled.
        Run-times are in seconds.
    }
    \label{tab:single-commodity}

    \centering

    \begin{threeparttable}
    \footnotesize

    \begin{tabular}{lrrrrrrrrrrrr}
    \toprule
     & \multicolumn{3}{l}{BB} & \multicolumn{3}{l}{FlowMIS} & \multicolumn{3}{l}{MIS} & \multicolumn{3}{l}{SNC} \\
     \cmidrule(lr){2-4}\cmidrule(lr){5-7}\cmidrule(lr){8-10}\cmidrule(lr){11-13}
     & No. iters & \#S & Time & No. iters & \#S & Time & No. iters & \#S & Time & No. iters & \#S & Time \\
    \midrule
    R04 & 35.5 & 18 & 94 & 15.3 & 18 & 15 & 11.5 & 18 & 10 & 12.1 & 18 & 11 \\
    R05 & 38.1 & 18 & 91 & 17.6 & 18 & 12 & 17.6 & 18 & 14 & 18.2 & 18 & 18 \\
    R06 & 14.7 & 18 & 33 & 9.9 & 18 & 8 & 9.0 & 18 & 7 & 9.1 & 18 & 6 \\
    R07 & 34.7 & 18 & 119 & 17.3 & 18 & 26 & 16.9 & 18 & 25 & 17.1 & 18 & 23 \\
    R08 & 26.6 & 18 & 68 & 11.3 & 18 & 2 & 10.0 & 18 & 4 & 10.0 & 18 & 4 \\
    R09 & 31.8 & 17 & 186 & 11.0 & 18 & 17 & 9.4 & 18 & 16 & 8.9 & 18 & 14 \\
    R10 & 24.2 & 18 & 89 & 13.4 & 18 & 17 & 12.8 & 18 & 5 & 12.6 & 18 & 5 \\
    \midrule
    Solved && 125 &&& \textbf{126} &&& \textbf{126} &&& \textbf{126} & \\
    Fastest && 9 &&& \textbf{48} &&& 34 &&& 35 & \\
    \multicolumn{2}{l}{Speed-up\tnote{1}} && &&& \textbf{8.7\texttimes} &&& 8.4\texttimes &&& 8.0\texttimes \\
    \bottomrule
    \end{tabular}
    \begin{tablenotes}
        \item[1]
        The speed-up is computed w.r.t. the BB formulation, for all other formulations. 
        It reports the average speed-up of the other formulations on instances both BB and the other formulation solve using the geometric mean~\citep{fleming_how_1986}.
    \end{tablenotes}
    \end{threeparttable}
\end{table}

\subsection{Results on multi-commodity instances}
\label{subsec:results_multicommodity}

In this section we evaluate the performance of the four different subproblem formulations presented in~\cref{subsec:feasibility_cuts}, the valid inequalities presented in~\cref{subsec:valid_inequalities}, and the added value of the metric inequality strengthening procedure.

\paragraph{Feasibility cuts}
\Cref{tab:nothing} shows that without the valid inequalities and metric inequality strengthening procedure, MIS and SNC outperform BB and FlowMIS.
In particular, SNC solves 117 instances within the two hour time limit, MIS 117, FlowMIS 116, while BB manages to solve only 83 instances.
Of the solved instances, SNC solves 51 the fastest, MIS 38, FlowMIS 28, while BB is never the fastest to find an optimal solution.
The speed-ups over BB reflect this difference in performance: while FlowMIS is 10.6 times faster than BB, MIS and SNC are 21.0 and 22.2 times faster to find a solution, respectively.
Instances in the groups R06 and R07 seem particularly challenging as no formulations manages to solve all instances in these groups within the time limit.

\Cref{tab:nothing} also presents results obtained using the large, deterministic equivalent (DEQ) formulation of~\eqref{eq:deq}.
Using this formulation, each instance is solved using four cores of a 2.45GHz AMD 7763 processor, 32GB of memory, and eight hours of run-time.
Despite the significant increase in resources, it is clear that the DEQ formulation is not competitive: it manages to solve only 81 instances, requiring on average nearly three times the run-time to solve an instance that the BB formulation does.

\begin{table}
    \caption{
        Average number of iterations, run-time, and number of solved multi-commodity instances (`\#S') for each subproblem formulation, by instance group, without the valid inequalities or the metric inequality strengthening procedure.
        Run-times are in seconds.
    }
    \label{tab:nothing}

    \centering

    \begin{threeparttable}
    \footnotesize
    \begin{tabular}{lrrrrrrrrrrrrrr}
    \toprule
     & \multicolumn{2}{l}{DEQ} & \multicolumn{3}{l}{BB} & \multicolumn{3}{l}{FlowMIS} & \multicolumn{3}{l}{MIS} & \multicolumn{3}{l}{SNC} \\
     \cmidrule(lr){2-3}\cmidrule(lr){4-6}\cmidrule(lr){7-9}\cmidrule(lr){10-12}\cmidrule(lr){13-15}
     & \#S & Time & No. iters & \#S & Time & No. iters & \#S & Time & No. iters & \#S & Time & No. iters & \#S & Time \\
    \midrule
    R04 & 16 & 3410 & 86.7 & 17 & 1018 & 33.3 & 18 & 133 & 29.1 & 18 & 104 & 25.2 & 18 & 74 \\
    R05 & 15 & 4358 & 62.3 & 18 & 409 & 17.2 & 18 & 70 & 14.2 & 18 & 36 & 11.6 & 18 & 24 \\
    R06 & 11 & 3654 & 106.9 & 12 & 583 & 83.8 & 15 & 593 & 71.0 & 16 & 422 & 72.2 & 16 & 498 \\
    R07 & 11 & 6117 & 172.0 & 6 & 916 & 71.8 & 12 & 968 & 56.2 & 11 & 473 & 65.8 & 12 & 569 \\
    R08 & 9 & 1711 & 179.2 & 12 & 1059 & 28.8 & 17 & 199 & 21.7 & 18 & 127 & 18.4 & 18 & 99 \\
    R09 & 12 & 6750 & 402.9 & 15 & 2334 & 21.5 & 18 & 188 & 14.8 & 18 & 105 & 13.6 & 18 & 122 \\
    R10 & 7 & 6862 & 711.7 & 3 & 3692 & 37.1 & 18 & 503 & 33.8 & 18 & 664 & 28.3 & 18 & 556 \\
    \midrule
    Solved & 81 &&& 83 &&& 116 &&& 117 &&& \textbf{118} & \\
    Fastest & 1 &&& 0 &&& 28 &&& 38 &&& \textbf{51} & \\
    \multicolumn{2}{l}{Speed-up\tnote{1}} & 0.7\texttimes &&& &&& 10.6\texttimes &&& 21.0\texttimes &&& \textbf{22.2\texttimes} \\
    \bottomrule
    \end{tabular}
    \begin{tablenotes}
        \item[1]
        The speed-up is computed w.r.t. the BB formulation, for all other formulations. 
        It reports the average speed-up of the other formulations on instances both BB and the other formulation solve using the geometric mean~\citep{fleming_how_1986}.
    \end{tablenotes}
    \end{threeparttable}
\end{table}

There does not appear to be an obvious pattern to the results presented in~\cref{tab:nothing}.
On the whole, SNC and MIS are faster for the ten-node instance groups R04 through R09.
FlowMIS is slightly faster on the larger twenty-node instance group R10.
Without additional strengthening, it is clear that either MIS or SNC should be the preferred formulation: these maximise the number of solved instances while obtaining the best speed-ups.

\paragraph{Valid inequalities}
\Cref{tab:only_vis} shows the results with the scenario creation technique of~\cref{subsec:valid_inequalities} enabled.
From now on we no longer include the DEQ in our tables, since the enhancements we will discuss do not apply to that formulation.
Comparing~\cref{tab:nothing} and~\cref{tab:only_vis}, it is clear that the scenario creation technique of~\cref{subsec:valid_inequalities} is very beneficial.
The number of iterations and average run-times both decrease substantially, in many cases by 5-10\texttimes, while several additional instances in the R06 group are now also solved within the time limit.
Thus, the valid inequalities are effective already in their own right.

Additionally, the performance characteristics of the different formulations change substantially with the additional structure the scenario provides.
\Cref{tab:only_vis} shows that FlowMIS and SNC formulations solves 119 instances with the valid inequalities in place, just one more than the 118 MIS solves, and two more than the 117 BB manages to solve.
FlowMIS is now clearly the faster formulation, solving 68 instances the fastest, against 33 for BB, while SNC and MIS are fastest on only 13 and 5 instances, respectively.
Although SNC and MIS both solve one or two more instances than BB, on the instances they both solve, BB is about as fast.
This is evidenced by the speed-ups over BB: for FlowMIS, this is 1.8 times, while MIS and SNC both obtain speed-ups of 1.0 times, and are thus on par with BB.

\begin{table}
    \caption{
        Average number of iterations, run-time, and number of solved multi-commodity instances (`\#S') for each subproblem formulation, by instance group, with the valid inequalities enabled.
        Run-times are in seconds.
    }
    \label{tab:only_vis}

    \centering

    \begin{threeparttable}
    \footnotesize

    \begin{tabular}{lrrrrrrrrrrrr}
    \toprule
     & \multicolumn{3}{l}{BB} & \multicolumn{3}{l}{FlowMIS} & \multicolumn{3}{l}{MIS} & \multicolumn{3}{l}{SNC} \\
     \cmidrule(lr){2-4}\cmidrule(lr){5-7}\cmidrule(lr){8-10}\cmidrule(lr){11-13}
     & No. iters & \#S & Time & No. iters & \#S & Time & No. iters & \#S & Time & No. iters & \#S & Time \\
    \midrule
    R04 & 19.7 & 18 & 120 & 20.6 & 18 & 142 & 19.6 & 18 & 111 & 18.9 & 18 & 112 \\
    R05 & 6.9 & 18 & 4 & 6.9 & 18 & 4 & 6.9 & 18 & 5 & 6.9 & 18 & 5 \\
    R06 & 63.9 & 17 & 911 & 55.8 & 17 & 918 & 51.0 & 16 & 681 & 59.5 & 17 & 778 \\
    R07 & 60.3 & 10 & 1117 & 55.9 & 12 & 520 & 58.0 & 12 & 734 & 52.0 & 12 & 568 \\
    R08 & 8.9 & 18 & 41 & 8.7 & 18 & 33 & 8.9 & 18 & 44 & 8.9 & 18 & 49 \\
    R09 & 4.1 & 18 & 7 & 4.1 & 18 & 7 & 4.1 & 18 & 8 & 4.1 & 18 & 8 \\
    R10 & 3.8 & 18 & 44 & 3.7 & 18 & 28 & 3.8 & 18 & 38 & 3.8 & 18 & 58 \\
    \midrule
    Solved && 117 &&& \textbf{119} &&& 118 &&& \textbf{119} & \\
    Fastest && 33 &&& \textbf{68} &&& 5 &&& 13 & \\
    \multicolumn{2}{l}{Speed-up\tnote{1}} && &&& \textbf{1.8\texttimes} &&& 1.0\texttimes &&& 1.0\texttimes \\
    \bottomrule
    \end{tabular}
    \begin{tablenotes}
        \item[1]
        The speed-up is computed w.r.t. the BB formulation, for all other formulations. 
        It reports the average speed-up of the other formulations on instances both BB and the other formulation solve using the geometric mean~\citep{fleming_how_1986}.
    \end{tablenotes}
    \end{threeparttable}
\end{table}

\paragraph{Metric inequalities}
Finally, we turn to the metric inequality strengthening procedure of~\citet{costa_benders_2009}.
\Cref{tab:everything} shows the results when these metric inequalities are also enabled, on top of the valid inequalities presented in~\cref{tab:only_vis}.
FlowMIS and MIS now solve 120 instances.
BB solves two instances less than in~\cref{tab:only_vis}, and is fastest on 31 instances.
FlowMIS remains the fastest formulation, solving 67 instances the fastest.
MIS and SNC are fastest for 12 and 10 instances, respectively.

Comparing run-times, there does not seem to be a meaningful performance improvement due to the strengthening procedure.
Although on the whole the number of iterations needed for convergence decreases somewhat for most formulations, the added cost of each iteration appears to cancel out most of the gains.

\begin{table}
    \caption{
        Average number of iterations, run-time, and number of solved multi-commodity instances (`\#S') for each subproblem formulation, by instance group, with the valid inequalities and metric inequality strengthening procedure enabled.
        Run-times are in seconds.
    }
    \label{tab:everything}

    \centering

    \begin{threeparttable}
    \footnotesize

    \begin{tabular}{lrrrrrrrrrrrr}
    \toprule
     & \multicolumn{3}{l}{BB} & \multicolumn{3}{l}{FlowMIS} & \multicolumn{3}{l}{MIS} & \multicolumn{3}{l}{SNC} \\
     \cmidrule(lr){2-4}\cmidrule(lr){5-7}\cmidrule(lr){8-10}\cmidrule(lr){11-13}
     & No. iters & \#S & Time & No. iters & \#S & Time & No. iters & \#S & Time & No. iters & \#S & Time \\
    \midrule
    R04 & 19.7 & 18 & 114 & 19.6 & 18 & 107 & 18.5 & 18 & 91 & 18.6 & 18 & 101 \\
    R05 & 6.9 & 18 & 4 & 6.9 & 18 & 4 & 6.9 & 18 & 5 & 6.9 & 18 & 4 \\
    R06 & 59.9 & 16 & 799 & 58.4 & 17 & 533 & 47.8 & 17 & 513 & 47.6 & 17 & 769 \\
    R07 & 61.1 & 9 & 880 & 57.4 & 13 & 660 & 57.5 & 13 & 727 & 53.6 & 12 & 790 \\
    R08 & 8.9 & 18 & 41 & 8.7 & 18 & 34 & 8.9 & 18 & 46 & 8.9 & 18 & 51 \\
    R09 & 4.1 & 18 & 7 & 4.1 & 18 & 7 & 4.1 & 18 & 9 & 4.1 & 18 & 8 \\
    R10 & 3.8 & 18 & 40 & 3.7 & 18 & 27 & 3.8 & 18 & 36 & 3.8 & 18 & 47 \\
    \midrule
    Solved && 115 &&& \textbf{120} &&& \textbf{120} &&& 119 & \\
    Fastest && 31 &&& \textbf{67} &&& 12 &&& 10 & \\
    \multicolumn{2}{l}{Speed-up\tnote{1}} && &&& \textbf{2.0\texttimes} &&& 1.0\texttimes &&& 1.0\texttimes \\
    \bottomrule
    \end{tabular}
    \begin{tablenotes}
        \item[1]
        The speed-up is computed w.r.t. the BB formulation, for all other formulations. 
        It reports the average speed-up of the other formulations on instances both BB and the other formulation solve using the geometric mean~\citep{fleming_how_1986}.
    \end{tablenotes}
    \end{threeparttable}
\end{table}

\vspace{1em}

Summarising, the additional structure of the valid inequalities allows FlowMIS to outperform the other formulations.
This suggests that (a) it is somewhat problematic to evaluate the solution approach's components in isolation (since the results of~\cref{tab:nothing} would have pointed us in a very different direction than the later tables), and (b) there is significant interaction between the various components.
Only when brought together do they yield their best results.

\section{Conclusion}
\label{sec:conclusion}

This paper studies a chance-constrained variant of the multicommodity capacitated network design problem in which the goal is to minimise the network construction costs, while ensuring that commodity flow can be routed through the network sufficiently often to meet uncertain commodity demand.
We propose an efficient solution approach based on Benders' decomposition.
We pay particular attention to the manner in which feasibility cuts are generated, and propose a new subproblem formulation that we term FlowMIS.
We show that FlowMIS finds minimum-capacity cuts in the single-commodity case, while other formulations typically do not.
We also add several valid inequalities to the master problem which are very beneficial in speeding up the overall convergence of the algorithm.
Numerical experiments on 126 benchmark instances demonstrate that FlowMIS formulation combined with the valid inequalities results in a competitive algorithm that outperforms the other subproblem formulations.
In particular, when all enhancements are enabled, FlowMIS solves 67 out of 120 solved instances the fastest, and achieves a speed-up of 2.0\texttimes~over a basic implementation, whereas the other subproblem formulations do not noticeably improve over this basic implementation.

The performance and simplicity of FlowMIS and our valid inequalities provide a solid foundation for further investigation of Benders' decomposition algorithms for network design problems.
In future work, one might investigate the theoretical properties of the FlowMIS feasibility cuts for multi-commodity problems.
Additionally, it will be interesting to see how our algorithm performs numerically on larger instances with even more scenarios.

\ACKNOWLEDGMENT{%
Niels A. Wouda and Evrim Ursavas gratefully acknowledge funding from the European Union’s Horizon 2020 research and innovation programme under grant agreement \#101022484.
}

\begin{APPENDICES}

\end{APPENDICES}

\bibliographystyle{informs2014trsc}
\bibliography{references}

\end{document}